\documentclass{amsart}
\usepackage{geometry, graphicx, layout, multicol, tikz, textcomp, wrapfig, float, xfrac, faktor, nicefrac, setspace, amsthm, epstopdf, xcolor, multirow, amssymb, amsmath, tensor}
\usepackage[utf8]{inputenc}
\usepackage[english]{babel}
\usepackage[utf8]{inputenc}
\usepackage[english]{babel}
\usepackage[pagewise]{lineno}
\theoremstyle{plain}
\usepackage{enumerate,bm}

\title{Groups with a Fixed Character Degree}
\author{Mark L. Lewis \and Brandon Martin}
\address{Department of Mathematical Sciences, Kent State University, Kent, OH 44242}
\email{lewis@math.kent.edu}
\email{bmarti52@kent.edu}

\begin{document}
\onehalfspacing
\newtheorem{lemma}{Lemma}
\newtheorem{theorem}{Theorem}
\newtheorem{sublemma}{Lemma}[lemma]
\theoremstyle{named}
\newtheorem*{theorem*}{Theorem}
\newtheorem*{namedtheorem}{Main Theorem}

\maketitle

\begin{abstract} This note is concerned with the question: for which positive integers $d,e$, with $d$ square-free and $\text{gcd}(d,d+e)=1$, does there exist a solvable (we will see, in our case, metabelian) group $G$, of order $|G|=d(d+e)$, such that $d$ is an irreducible character degree of $G$?  The main theorem is an arithmetic characterization in terms of the prime divisors of $d(d+e)$.  This work is motivated by previous investigations into finite groups $G$ that have an irreducible character of comparatively large degree.
\end{abstract}

\section{Introduction} 
We will let all groups be finite throughout this note.  Let $G$ be a group of order $d(d+e)$ where $e>1$ is an integer and $d\in\text{cd}(G)$, where $\text{cd}(G)=\{\chi(1)\text{ $|$ }\chi\in\text{Irr}(G)\}$.  Lewis \cite{Lew} proved that $|G|\leq e^4-e^3$ for groups with a nontrivial abelian normal subgroup.  Following this, Nguyen, Lewis, and Schaeffer Fry \cite{LHS} proved that if $d$ is the degree of a complex irreducible character of $G$ where $|G|=d(d+e)$ for some integer $e>1$, then $|G|\leq e^4-e^3$.  Their proof uses the classification of finite simple groups.  Also, by earlier work by Isaacs \cite{Is}, the bound of $e^4-e^3$ is the best possible. \par 

The above results assume that $d\in\text{cd}(G)$.  We now consider this problem from the opposite perspective, as suggested by the authors of \cite{LHS}.  In particular, we will no longer assume that $d\in\text{cd}(G)$, and we look to characterize conditions on $G$ which result in $d$ being an irreducible character degree of $G$. \par 

Let $G$ be a group of order $d(d+e)$, where $d,e>1$.  Snyder \cite{Sny} classified all such groups with $e=2$ and $e=3$ which has $d$ as a complex irreducible character degree.
Durfee and Jensen \cite{DJ} classified all such groups with $4\leq e\leq 6$, as well as all irreducible character degrees, $d$, which can occur with $e=7$. Following this, Sambale \cite{Sam} extended this classification of $d$ for all values when $e\leq 11$.  Much of the previous work has relied on computer algebra systems, with obvious limitations as $e$ gets larger.   \par
We seek to answer this question for all solvable groups of order $d(d+e)$ by finding an arithmetic characterization in terms of the prime divisors of $d$ and $d+e$ for when $d\in\text{cd}(G)$. In this note, we only consider the case where $d$ and $d+e$ are coprime and $d$ is square-free.

\begin{namedtheorem}
There exists a metabelian group $G$, such that \newline $|G|=d_1\cdots d_mp_1^{a_1}\cdots p_n^{a_n}$, where the set $\{d_1,\dots,d_m,p_1,\dots,p_n\}$ are distinct primes and $a_i\in\mathbb{N}$, such that $d=d_1\cdots d_m\in\text{cd}(G)$, if and only if there exist $\epsilon_1,\dots, \epsilon_k\in\mathbb{N}$, and $b_1,\dots,b_k\in\mathbb{N}$, where $k\leq m$, such that $$q_1^{\epsilon_1}\equiv 1({\rm mod}\text{ }b_1),\dots, q_k^{\epsilon_k}\equiv 1({\rm mod}\text{ }b_k),$$ where $b_1\cdots b_k=d$ and primes $q_1,\dots,q_k$ such that $q_1^{\epsilon_1}\dots q_k^{\epsilon_k}\mid p_1^{a_1}\dots p_n^{a_n}$.
\end{namedtheorem} 
We make a few remarks to clarify the notation of the Main Theorem:
\begin{enumerate}
    \item Each $q_r$ need not be unique for any $1\leq r\leq k$.  If $q_\alpha^{e_\alpha}=q_\beta^{\epsilon_\beta}$ for some $\alpha,\beta\in\{1,\dots, k\}$, we only require that $q_\alpha^{e_\alpha}q_\beta^{\epsilon_\beta}\mid p_i^{a_i}$ where $i=\alpha=\beta$.
    \item We do not require $b_l$ to be prime for any $1\leq l\leq k$.
\end{enumerate}

The work in this paper was completed by the second author (Ph.D candidate) under the supervision of the first author at Kent State University.  The contents of this paper may appear as part of the second author's Ph.D dissertation.  The second author also thanks Dr. Stephen Gagola for his invaluable comments and suggestions.   

\section{Background}
We will need to make use of Frobenius groups throughout the course of this paper.  From \cite{CTFG} and \cite{FGT}, chapters $7$ and $3$, respectively, we collect the following known facts:  Let $G$ be a Frobenius group with Frobenius complement $H$.  Then there exists $K\triangleleft G$ with $HK=G$ and $H\cap K=1$.  Here, $K$ is called the Frobenius kernel of $G$.  We will continue to use this notation for $K$ and $H$ in this section.  We also know that Frobenius complements and Frobenius kernels necessarily have coprime orders. Indeed, $|K|\equiv 1(\text{mod }|H|)$.  \par 
The irreducible characters of Frobenius groups have been extensively studied, and they are of two types: ones that are induced from non-trivial irreducible characters of $K$, and others whose kernels contain $K$.  Suppose $\chi\in\text{Irr}(G)$ with $K\not\subseteq\text{ker}\chi$.  Then we have $\chi=\phi^G$ for some $\phi\in\text{Irr}(K)$.  It follows that $\phi^G(1)=|H|\phi(1)$.  So, if $\phi$ is a linear character, then $|H|=\phi^G(1)$. In particular, $|H|\in\text{cd}(G)$.  \par
By Thompson's classical result, we know that $K$ is nilpotent, thus $K'\lneq K$ where $K'$ is the derived subgroup of $K$.  Therefore, $K$ must have at least one linear character that is not the principal character. Let $\phi\in\text{Irr}(K)$ be such a character.  We know that $\phi$ lies in an orbit of size $|H|$ and $\phi^G\in\text{Irr}(G)$.  Hence $\phi^G(1)=|H|$.  In particular, we will always be guaranteed such a $\phi$ so that $|H|\in\text{cd}(G)$. \par    
From \cite{Mich}, we will state the It$\hat{\text{o}}$-Michler theorem, which will be used throughout the proof of the main theorem.  However, note that only It$\hat{\text{o}}$'s contribution will be needed for the proof of our result, and so our result does not rely on the classification of finite simple groups.
\begin{theorem*}[It$\hat{\text{o}}$-Michler]\label{Ito}
Let $G$ be a finite group and $p$ a prime.  Then $p\nmid\chi(1)$ for all $\chi\in\text{Irr}(G)$, if and only if $G$ has a normal, abelian Sylow $p$-subgroup.
\end{theorem*}

\section{Main Result}
We are now ready to proceed with the proof of the main theorem.
\begin{proof}[Proof of the Main Theorem]  
        First, we consider the backward direction. Suppose there exist $\epsilon_1,\dots,\epsilon_k\in\mathbb{N}$, and $b_1,\dots,b_k\in\mathbb{N}$ such that $$q_1^{\epsilon_1}\equiv 1({\rm mod}\text{ }b_1),\dots, q_k^{\epsilon_k}\equiv 1({\rm mod}\text{ }b_k),$$ where each $q_1,\dots,q_k$ is (not necessarily distinct) prime such that $q_1^{\epsilon_1}\dots q_k^{\epsilon_k}\mid p_1^{a_1}\dots p_n^{a_n}$, where each $p_i$ is (a distinct) prime.  Let $b_1\cdots b_k=d_1\cdots d_m=d$, where each $d_j$ is prime, and so we will let $1\leq k\leq m$.  
        
        Given the above sequence of congruences, it is straightforward to construct a group $G$ in the following way: consider the field $\mathbb{F}_{{q_i}^{\epsilon_i}}$. Consider the subgroup $B_i\leq \mathbb{F}_{{q_i}^{\epsilon_i}}^\times$ where $|B_i|=b_i$.  We have that $B_i$ acts naturally on the additive group $(\mathbb{F}_{{q_i}^{\epsilon_i}}, +)$ by multiplication in $\mathbb{F}_{{q_i}^{\epsilon_i}}$, and this is a Frobenius action with a cyclic Frobenius complement $B_i$ and Frobenius kernel $(\mathbb{F}_{{q_i}^{\epsilon_i}}, +)$.  Since $\mathbb{F}_{{q_i}^{\epsilon_i}}$ has prime characteristic, $(\mathbb{F}_{{q_i}^{\epsilon_i}}, +)$ is forced to be elementary abelian.  Thus each of our given congruences results in a Frobenius group with Frobenius complement of order $b_i$.  Taking the direct product of these Frobenius groups, as well as, if necessary, an abelian group of order $\frac{p_1^{a_1}\dots p_n^{a_n}}{q_1^{\epsilon_1}\dots q_k^{\epsilon_k}}$, we see that $b_1\cdots b_k=d\in\text{cd}(G)$, as wanted, by Theorem 4.21 from \cite{CTFG}.  Moreover, since $G$ is constructed as a direct product of cyclic groups acting on elementary abelian groups, we have that $G$ is metabelian.\par

        For the forward direction, we will now let $G$ be a metabelian (thus, solvable) group of order $d_1\cdots d_mp_1^{a_1}\cdots p_n^{a_n}$, where the set of primes \newline $\{d_1,\dots,d_m,p_1,\dots,p_n\}$ are all distinct, and $a_i\geq 1$ for all $i$.  Suppose $d=d_1\cdots d_m\in\text{cd}(G)$. \par

        We will proceed by induction on $m$.  First, consider the case when $m=1$, that is, when $d=d_1$ is prime.  We will proceed by contradiction.  Suppose $p_i^{\epsilon_i}\not\equiv 1(\text{mod }d_1)$ for each $i$ and for all $1\leq \epsilon_i\leq a_i$.  Fix $D\in\text{Syl}_d(G)$, a Sylow $d$-subgroup of $G$, and also fix $i$.  Since $G$ is solvable, we can choose $\hat{P_i}\in\text{Hall}_{dp_i}(G)$, a Hall $dp_i$-subgroup of $G$, such that $D\subseteq \hat{P_i}$.  By our assumption, we have that $D\trianglelefteq\hat{P_i}$. Let $P_i\in\text{Syl}_{p_i}(\hat{P_i})$, a Sylow $p_i$-subgroup of $\hat{P_i}$.  Then, $P_i\subseteq N_G(D)$.  Since $i$ was arbitrary, $P_i\subseteq N_G(D)$ for all $i$, and so $D\trianglelefteq G$, a contradiction to the It$\hat{\text{o}}$-Michler theorem.\par 
        
        Let $\chi\in\text{Irr}(G)$ such that $\chi(1)=d=d_1\cdots d_m$.  Let $N$ be a minimal normal subgroup of $G$.  Since $G$ is solvable, $N$ must be elementary abelian. If $d_j\mid |N|$ for any $j$, then $d\not\in\text{cd}(G)$, by \cite{Mich}.  Thus, $N$ is necessarily an elementary abelian $p$-group for some prime $p_v\mid p_1^{a_1}\cdots p_n^{a_n}$., where $v\in\{1,\dots, n\}$.  \par

        By Clifford's Theorem, $\chi_N=e\Sigma_{x=1}^t\theta_x$, where $\theta_1,\dots,\theta_t$ are conjugate irreducible characters of $N$ and $e=[\chi_N,\theta]$, for some  fixed irreducible constituent $\theta\in\text{Irr}(N)$ of $\chi_N$.  By Clifford, there exists a unique $\beta\in\text{Irr}(I_G(\theta)|\theta)$ such that $\beta^G=\chi$ and $\beta$ is the Clifford correspondent of $\chi$ in $I_G(\theta)$.  Since $\beta$ is the Clifford correspondent of $\chi$, $\beta(1)$ must be the full $d-$part of $|I_G(\theta)|$, and so we may apply the inductive hypothesis to $I_G(\theta)$.  Thus, suppose $|I_G(\theta)|=\delta\pi$, where $\delta$ is a divisor of $d_1\cdots d_m$, and $\pi$ is a divisor of $p_1^{a_1}\cdots p_n^{a_n}$.  Note that we've shown that $\beta(1)=\delta$.  Then, by induction, there exist $\epsilon_1,\dots, \epsilon _c\in\mathbb{N}$ and $b_1,\dots, b_c\in\mathbb{N}$, where $c\leq m$, such that 
         \begin{equation}
        q_1^{\epsilon_1}\equiv 1({\rm mod}\text{ }b_1), \dots, q_c^{\epsilon_c}\equiv 1({\rm mod}\text{ }b_c),
        \end{equation} where $b_1\cdots b_c=\delta$ and primes $q_1,\dots,q_c$ such that $q_1^{\epsilon_1}\cdots q_c^{\epsilon_c}\mid \pi$. \par Consider the set $\Delta=\{d_j\text{ : }d_j\mid |G:I_G(\theta)|\}$.  Let $D_j\in\text{Syl}_{d_j}(G)$, a Sylow $d_j$-subgroup of $G$, and let $D\in\text{Hall}_\Delta(G)$, a Hall $\Delta$-subgroup of $G$.  Note that for each $j$, $D_j$ is cyclic of prime order $d_j$.  Since $C_D(N)=I_G(\theta)\cap D=1$, we have that $C_D(N)=1$ and so $D$ acts faithfully on $N$.  \par

        If $C_N(D)\neq 1$ then $C_N(D)\trianglelefteq G$ centralizes $N$.  Since $N$ is a minimal normal subgroup of $G$, our assumption forces $C_N(D)=N$, which is a contradiction as $|N|$ and $|D|$ are coprime.  Therefore, $C_N(D)=1$ and so $N$ is not central.\par 
        Now, we may write $N=N_1\times\cdots\times N_l$, where each $N_i$ is an elementary abelian $p_v$-group, and each $N_i$ is irreducible under the action of $D$, by Maschke's theorem.  We know that each $D_j$, where $d_j\in\Delta$, acts Frobeniusly on at least one factor of $N$. If not, then such a $D_j$ will centralize each $N_i$, a contradiction to $D$ acting faithfully on $N$. Suppose, after reordering, that we have $$\Delta=\{d_1,\dots,d_w\}.$$ Let 
        $$\Delta_1=\{j\in\{1,\dots, w\} \mid D_j \text{ acts Frobeniusly on $N_1$}\}.$$  
        This gives a Frobenius group with Frobenius complement $\prod_{j\in\Delta_1}D_j$ and Frobenius kernel $N_1$.  Next, we will let $$\Delta_2=\{j\in\{1,\cdots,w\}\setminus\Delta_1\mid D_j\text{ acts Frobeniusly on $N_2$}\}.$$ 
        This gives another Frobenius group with Frobenius complement $\prod_{j\in\Delta_2}D_j$ and Frobenius kernel $N_2$.  We continue in this way and let $$\Delta_s=\{j\in\{1,\dots,w\}\setminus\bigcup_{z=1}^{s-1}\Delta_z\mid D_j\text{ acts Frobeniusly on $N_s$}\}.$$  Therefore, we get the following sequence of congruences: 
        \begin{equation} |N_1|\equiv 1(\text{mod }|\prod_{j\in\Delta_1}D_j|),\dots, |N_s|\equiv 1(\text{mod }|\prod_{j\in\Delta_s}D_j|),\dots,\end{equation} where $|\prod_{j\in\Delta_1}D_j|\cdots|\prod_{j\in\Delta_s}D_j|\cdots =|D|$ and each $N_i$ is a $p_v$-group such that $|N_1|\cdots |N_s|\cdots$ divides $|N|$. Note that there may be factors of $N$ that do not appear in our sequence of congruences.  We know that there must exist some $D_j$, such that $d_j\in\Delta$, where $D_j$ acts Frobeniusly on such a factor of $N$. However, we do not need to consider this resulting congruence, because that particular Sylow $d_j$-subgroup has already been accounted for due to its action on a different factor of $N$.  By $(1)$ and $(2)$ above, we obtain the desired result.\par 
        Suppose $G=I_G(\theta)$.  Then $\ker(\theta)$ is G-invariant and so $\ker(\theta)\trianglelefteq G$ and $\ker(\theta)\leq N$. Since $N$ is minimal normal, we must have $\ker(\theta)=1$ or $N$.  Since, by Clifford, $\theta\neq 1_N$, we know $\ker(\theta)\neq N$ and so $\ker(\theta)=1$.  Thus, $\theta$ is faithful and since $N$ is elementary abelian, we are forced to have $N$ be cyclic.  Therefore,  $N$ must have prime order.  Now, since $\theta$ is $G$-invariant, every conjugation automorphism fixes $\theta$.  Since $\theta$ is faithful and $N$ is a cyclic group of prime order, the only automorphism that fixes $\theta$ is the identity.  Thus, every $g\in G$ centralizes $N$ and so $N\leq Z(G)$.  Since $N\leq Z(G)$, Schur's lemma implies $\chi_N=\chi(1)\lambda$ for some unique $\lambda\in\text{Irr}(N)$, and, in particular, here $\lambda=\theta$.  Thus, $d=\chi(1)\in\text{cd}(G/N)$, and so we apply the inductive hypothesis to $G/N$ to obtain the result.   
\end{proof}

\end{document}